\documentclass[11pt,a4paper,twoside]{amsart}
\usepackage{amsmath}
\usepackage{amssymb}
\usepackage{latexsym}

\newcommand{\co}{\colon\thinspace}    
\newcommand{\fnote}[1]{\footnote{\small sharp1}}
\newcommand{\inv}{^{-1}}              

\newcommand{\N}{{\mathbb N}}
\newcommand{\Z}{{\mathbb Z}}
\newcommand{\R}{{\mathbb R}}
\newcommand{\Q}{{\mathbb Q}}

\newcommand{\T}{{\mathbb T}}

\newcommand{\spt}{\mbox{supp}}
\newcommand{\Azero}{\mathcal{A}_0}

\newtheorem{theorem}{Theorem}[section]
\newtheorem{proposition}[theorem]{Proposition}
\newtheorem{corollary}[theorem]{Corollary}
\newtheorem{definition}[theorem]{Definition}
\newtheorem{lemma}[theorem]{Lemma}

\newtheorem{problem}[theorem]{Problem} 

\title{Differentiability of Mather's $\beta$-function vs Ma\~n\'e's conjecture}

\author{Daniel Massart}
\date{\today}
\begin{document}

\begin{abstract}
We prove that if a time-periodic Tonelli Lagrangian on a closed manifold $M$ satisfies a strong version of the Differentiability Problem for Mather's $\beta$-function, then the Legendre transforms of rational homology classes are dense in the first cohomology of $M$, which is a first step towards Ma\~n\'e's conjecture.
\end{abstract}
\maketitle

\section{Introduction}
\subsection{Tonelli Lagrangians}
Let $M$ be a compact, connected manifold without boundary. For the sake of brevity we shall refer to such manifolds as closed.   A Tonelli Lagrangian on $M$ is a $C^{2}$ function on $TM\times \T$, where $\T$ is the  circle $\R / \Z$, satisfying the following conditions : 
\begin{enumerate}
\item
for every $(x,t) \in M \times  \T$, the function $v \mapsto L(x,v,t)$ 
is superlinear
\item
for every $(x,v,t) \in TM \times  \T$, 
the bilinear form $\partial^2 L(x,v,t) / \partial v^2$ is positive definite 
\item
the local flow $\Phi_t$  defined on $TM\times \T$  
by the Euler-Lagrange equation for extremals of the action of curves 
is complete.
\end{enumerate}

  A good  example to keep in mind is the sum of a Riemann  metric, viewed as a quadratic function on $TM$, and a time-periodic potential (a function on $M\times \T$).   See \cite{Mather91, Fathi}  for more background and references. When the Lagrangian does not depend on $t \in \T$, it is called autonomous and we omit the factor $\T$. 

The following classical way  to obtain invariant subsets of the Euler-Lagrange flow was introduced by Mather in \cite{Mather91}.  Define $\mathcal{M}_{inv}$
to be the set of $\Phi_t$-invariant, compactly supported, Borel probability measures on $TM\times \T$.
Mather showed that the function (called action of the Lagrangian on measures)
	\[
	\begin{array}{rcl}
\mathcal{M}_{inv} & \longrightarrow & \R \\
\mu & \longmapsto & \int_{TM\times \T}	L d\mu
\end{array}
\]
is well defined and has a minimum.  A measure achieving this minimum is called $L$-minimizing. The union of the supports of all minimizing measures is called Mather set, and denoted $\mathcal{M}(L)$. The following classical trick gives us more invariant sets from the same construction.  If $\omega$ is a closed one-form on  $M$, then $L-\omega$ is again a Tonelli  Lagrangian, and it has the same Euler-lagrange flow as $L$. Besides, by \cite{Mather91},  if $\mu \in \mathcal{M}_{inv}$, the integral $\int_{TM\times \T}	\omega d\mu$ only depends on the cohomology class of $\omega$.  A measure achieving the minimum of 
	\begin{equation}\label{definition alpha}
	\begin{array}{rcl}
\mathcal{M}_{inv} & \longrightarrow & \R \\
\mu & \longmapsto & \int_{TM\times \T}	(L-\omega) d\mu
\end{array}
\end{equation}
where $\omega$ is any closed 1-form in the  cohomology class $c$, is called $(L,c)$-minimizing. Thus, for each cohomology class $c \in H^1(M,\R)$, we find a Mather set $\mathcal{M}(L,c)$.
\subsection{Ma\~n\'e's problem}
The upside of this measure-theoretic construction is that it readily yields an existence theorem. The downside is that we don't know what the minimizing measures look like. One of the most basic questions one may ask is, \emph{ if we choose the Lagrangian randomly, what are the minimizing measures ? are they supported on fixed points or periodic orbits ? if not, how big are their supports }? 
 
 In \cite{Mane95, Mane96}, R. Ma\~n\'e proposed the following problem (sometimes stated as Ma\~n\'e's conjecture, although Ma\~n\'e's conjecture is a stronger statement published in a later paper \cite{Mane97}): 
\begin{problem}\label{probleme Mane}
Is it true that given a Tonelli Lagrangian $L$ on a manifold $M$, there exists an residual subset $\mathcal{O}(L)$ of $C^{\infty}(M\times \T)$, such that for any 
$f \in  \mathcal{O}(L)$, there exists an open dense set $U(L,f)$ of $H^1(M,\R)$, such that for any $c \in U(L,f)$, there exists only one $(L+f,c)$-minimizing measure, and it is supported by a periodic orbit ?
\end{problem}
The case when $M$ is the circle is treated in \cite{Osuna}, with some arguments provided by \cite{Mather02}. The analogous problem in codimension one Aubry-Mather theory is dealt with in \cite{Bessi-Massart}. However, the full statement of Problem \ref{probleme Mane} seems way out of reach for the time being, so it makes sense to look for similar, but simpler, problems. 

First, we leave aside the question of the openness, and focus on finding a dense set of cohomology classes. 
Second, we relax the requirement on the support of the minimizing measure, by replacing it with some condition on its homology class. Let us explain what the homology class of a measure is. 
For any $\mu \in  \mathcal{M}_{inv}$, the homology class  $\left[\mu\right]$ is the unique  $h \in H_1 (M,\R)$ such that  
	\[
\langle \left[\omega \right],h \rangle = \int_{TM\times \T}	\omega d\mu 
\]
for any closed one-form $\omega$ on $M$. This is well defined because if $\mu \in \mathcal{M}_{inv}$, the integral $\int_{TM\times \T}	\omega d\mu$ only depends on the cohomology class of $\omega$. 

Now let us explain the notion of irrationality of a homology class. The torsion-free part of $H_1 (M,\Z)$ embeds as a lattice $\Gamma$ in $H_1 (M,\R)$. A class $h \in H_1 (M,\R)$ is called integer if it lies in $\Gamma$, and rational if $nh \in \Gamma$ for some $n \in \Z$. 
A subspace of $H_1 (M,\R)$ is called integer if it is generated by integer classes.

The quotient $H_1 (M,\R)/ \Gamma$ is a torus $\T^{b}$, where $b$ is the first Betti number of $M$. 
For $h$ in $H_1 (M,\R)$, the image of $\Z h$ in $\T^{b}$ is a subgroup of $\T^{b}$, 
hence its closure $\mathcal{T}(h)$ is a finite union of tori of equal dimension. This dimension is called the irrationality 
$I(h)$ of $h$ (it is denoted $I_{\Z}(h)$ in \cite{vmb2}). 
It is zero if  $h$ is rational. We say a class $h$ is completely irrational 
if its irrationality is maximal, i.e. equals $b$. 
In the same way, if $v$ is a vector of $\R^n$, we call irrationality of $v$ 
the dimension of the image of $\Z v$ in $\R^n / \Z^n$. 
Note that the irrationality of $h$ equals that of $nh$ for $n\in \Z$, 
$n \neq 0$ since the quotient of $\mathcal{T}(h)$ 
by $\mathcal{T}(nh)$ is a group of cardinality $n$. 

As an example let us look at the homology class of a measure supported on a periodic orbit. 
Assume $\gamma \co \left[0,T\right] \longrightarrow M$, with $T \in \N$,  is a $C^2$ closed curve such that $(\gamma, \dot\gamma, t)$ is a periodic orbit of $\Phi_t$. Then the probability measure
$\mu_{\gamma}$ on $TM\times \T$ such that, for any continuous function $f$ on $TM\times \T$, 
$$
\int_{TM} f d \mu_{\gamma} = \frac{1}{T}\int_{0}^{T} f(\gamma (t), \dot\gamma (t),t ) dt
$$
is $\Phi_t$-invariant and supported on $(\gamma, \dot\gamma, t)$, so it belongs in $\mathcal{M}_{inv}$. 
Its homology class is $ \frac{1}{T}\left[\gamma \right]$. Therefore $\left[\mu_{\gamma}\right]$ is rational. Conversely, one may ask whether an invariant measure with a rational  homology class is supported on a periodic orbit. It is true when $M=\T$, or when $\dim M= 2$ and $L$ is autonomous (see \cite{CMP}), however, it is false in general.
Still, a good starting point for Problem \ref{probleme Mane} is the following
\begin{problem}\label{probleme faible}
When is it true that given a Tonelli Lagrangian $L$ on a manifold $M$,  there exists a dense set $U(L)$ of $H^1(M,\R)$, such that for any $c \in U(L)$,  there exists a $(L,c)$-minimizing measure with a rational homology class ?
\end{problem}
\subsection{The Differentiability Problem}

In this paper we mean to explain the relationship between Problem \ref{probleme faible} and another important question in the field, the differentiability problem for Mather's $\beta$-function. First we need to define some objects. Mather's $\alpha$-function is the opposite of the minimum in Equation (\ref{definition alpha}), that is, 
\[
	\begin{array}{rcl}

\alpha \co H^1 (M,\R) & \longrightarrow & \R \\
c & \longmapsto & 
-\min \left\{\int_{TM\times \T}	(L-\omega)d\mu \co \mu \in  \mathcal{M}_{inv},\  \left[\omega\right]=c\right\}.
\end{array}
\]
This is a convex and superlinear function (see \cite{Mather91}), so it has a convex dual 
$$
\beta \co H_1(M,\R) \longrightarrow \R,
$$ 
defined by 
$$
\forall h \in H_1 (M,\R),\ \beta (h) = \sup_{c \in H^1 (M,\R)}\left( \langle c,h \rangle -\alpha (c) \right) .
$$
Actually Mather (\cite{Mather91}) proved that 
$$
\forall h \in H_1 (M,\R),\ \beta (h) = \min  \left\{\int_{TM\times \T}	 Ld\mu \co \mu \in  \mathcal{M}_{inv},\  \left[ \mu \right]= h \right\}.
$$
By the superlinarity of $\alpha$ and $\beta$, the suprema in the definitions of  $\alpha$ and $\beta$ are actually maxima. 
In particular $\min \alpha = -\beta(0)$, and for all $c \in H^1(M,\R), h \in H_1(M,\R)$ we have the Fenchel inequality :
$$
\alpha(c)+\beta(h) \geq \langle c,h \rangle .
$$
For any $h \in H_1(M,\R)$, we call Legendre transform of $h$ with respect to $\beta$, the set
$$
\partial \beta (h) := \{ c \in H^1(M,\R) \co \alpha(c)+\beta(h) =  \langle c,h \rangle \}.
$$
Likewise, for any $c \in H^1(M,\R)$, we call Legendre transform of $c$ with respect to $\alpha$, the set
$$
\partial \alpha (c) := \{ h \in H_1(M,\R) \co \alpha(c)+\beta(h) =  \langle c,h \rangle \}.
$$
The functions $\alpha$ and $\beta$ are sometimes called effective Hamiltonian and Lagrangian, respectively. They resemble the Hamiltonian, or Lagrangian, in that they are convex and superlinear; on the other hand, they need not be strictly convex, nor smooth, that is, the sets $\partial \beta (h)$
and $\partial \alpha (c)$ need not have cardinality one. By convex duality, it is equivalent to study the differentiability of $\beta$, and the strict convexity of $\alpha$. 

When $M=\T$, the following theorem says everything about the differentiability of $\beta$ :  
\begin{theorem}[\cite{Mather90, Bangert94}]\label{Mather}
If $M=\T$ then $\beta$ is differentiable at every irrational homology class. 
It is differentiable at a rational homology class if and only if periodic orbits in this class fill up $\T$. 
\end{theorem}
If we want to extend this theorem to higher dimensional manifolds, it seems natural to proceed as follows. 
A convex function has a tangent cone at every point. We say that $\beta$ is differentiable in $k$ directions at $h$ if the tangent cone to $\beta$ at $h$ 
contains a linear space of dimension $k$. 
We are thus led to ask \textit{whether $\beta$ is always differentiable in at least $k$ directions at a $k$-irrational homology class}.  This will henceforth be referred to as the Differentiability Problem.
Mather conjectures the answer is yes for $C^{\infty}$ Lagrangians. 
The answer to the Differentiability Problem is yes   for all $C^{2}$ Lagrangians when $M= \T$ by Theorem \ref{Mather}. It cannot be yes in general by \cite{BIK}. 

\subsection{A stronger version of the Differentiability Problem}
\begin{definition}\label{def V_h}
Let \begin{itemize}
  \item $M$ be a closed manifold
  \item $L$ be a Tonelli Lagrangian on $TM \times \T$
  \item $h$ be a homology class in $H_1(M,\R)$.
\end{itemize}
We define 
$$
\tilde{V}_h := \{ \lambda \left(c-c', \alpha(c)- \alpha(c')\right) \co c,c' \in \partial \beta (h), \  \lambda \in \R \}.
$$
\end{definition}
In plain language, $\tilde{V}_h $ is the underlying vector space to the affine space generated by the subset $\{(c, \alpha(c)) \co c \in \partial \beta (h) \}$ of $H^1(M,\R) \times \R$.
We define 
$$
V_h := \{ \lambda \left(c-c \right) \co c,c' \in \partial \beta (h), \  \lambda \in \R \},
$$
that is, $V_h$ is the canonical projection of $\tilde{V}_h $ to $H^1(M,\R) $.

Recall that a vector subspace of $H^1 (M,\R)\times \R = H^1 (M\times \T, \R)$ is integer if it is generated by integer cohomology classes. A slight modification of  \cite{vmb2}, Proposition 20 (see Proposition \ref{prop vmb2, modifiee}) says  that \emph{if for every $h \in H_1 (M,\R)$, $\tilde{V}_h$ is integer, then the answer to the Differentiability Problem is affirmative, that is, $\beta$ is always differentiable in at least $k$ directions at a $k$-irrational homology class}.
\subsection{An even stronger version of the Differentiability Problem}\label{stronger stronger}
In \cite{ijm} we proposed the following approach to  the Differentiability Problem, based on the notion of Aubry set, as defined by Fathi in \cite{Fathi} (see subsection \ref{Aubry}). For the time being all we need to know is that given a Tonelli Lagrangian $L \co TM \times \T \longrightarrow \R$ and an homology class $h$, the Aubry set $\tilde{\mathcal{A}}(L,h)$ is a compact subset of $TM \times \T$ which is invariant under the Euler-Lagrange flow of $L$, and projects injectively to $M\times \T$. The projection of $\tilde{\mathcal{A}}(L,h)$ to $M\times \T$ is called the projected Aubry set, and denoted $\mathcal{A}(L,h)$. Another  important object in weak KAM theory is the so-called quotient Aubry set, originally defined in \cite{Mather02} (see also \cite{FFR}, and \cite{Bernard_Conley}). See subsubsection \ref{Aubry quotient} for the precise definition. 

The idea of \cite{ijm} is to relate the differentiability of $\beta$ with the topology of the complement in $M\times \T$ of the projected Aubry set. 

Let $\tilde{E}_h$ be   the set of $(c,\tau) \in H^1 (M\times \T,\R)=H^1 (M,\R)\times H^1 ( \T,\R)$ such that there exists  a smooth closed one-form $\omega$ on $M\times \T$ with $[\omega] = (c,\tau)$ and $\spt (\omega)\cap \tilde{\mathcal{A}}(L,h) = \emptyset$. Let $E_h$ be the canonical projection of $\tilde{E}_h$ to $H^1 (M,\R)$. By a slight modification (see subsubsection \ref{E_h et E_c}) of   Theorem 14 of  \cite{soussol} (\cite{ijm} in the autonomous case) we have,  for every $h \in H_1 (M,\R)$, $E_h \subset V_h$.

What this has to do with the Differentiability Problem is summed up in  Corollary \ref{differentiabilite et E_h = V_h} which  says that if, for every $h \in H_1 (M,\R)$, we have $E_h=V_h$, then we have  an affirmative answer to  the Differentiability Problem.
\subsection{Now at last we can state our main theorem}
The meaning of our main result is, roughly speaking, that if we have an affirmative answer to the strong version of the Differentiability Problem stated in Subsection \ref{stronger stronger}, plus some non-degenaracy conditions which are satisfied by a generic Lagrangian thanks to \cite{BC}, or when the dimension of $M$ is small by \cite{FFR}, then we  have an affirmative answer to Problem \ref{probleme faible}, which is  a weak version of Problem \ref{probleme Mane}, which is itself a weak version of Ma\~n\'e's conjecture. All of which is to say that Ma\~n\'e's conjecture is a pretty bold statement. In \cite{FR}, another approach to Ma\~n\'e's conjecture is taken, by relaxing the regularity condition : instead of looking for a residual subset of $C^{\infty}(M)$, the authors look for a residual subset of $C^{2}(M)$. 

Before stating our result, let us point out that given a cohomology class $c$, it is equivalent to say that there exists a $(L,c)$-minimizing measure with a rational homology class, and to say that there exists a rational homology class $h$ such that $c \in  \partial \beta (h)$. With that in mind, the following theorem partially answers Problem \ref{probleme faible}.
\begin{theorem}\label{main}
Let  $M$ be a closed manifold, and let $L$ be a Tonelli Lagrangian on $TM \times \T$. Assume that  for every $h \in H_1 (M,\R)$, $V_h=  E_h$, and the quotient Aubry set $A_h$ has Hausdorff one-dimensional measure zero.
 Then the set of cohomology classes 
 $$
 \bigcup \partial \beta (h) \co h \in H_1 (M,\R), h \mbox{ rational}
 $$
 is dense in $H^1(M,\R)$.
\end{theorem}
This theorem is proved, in the context of codimension one Aubry-Mather theory, as Proposition 2.2 of \cite{Bessi-Massart}. In the codimension one theory, the non-degeneracy hypothesis-that is, the quotient Aubry set being small- is not necessary, and $\alpha$ is $C^1$, which makes the proof a little easier. 
\section{Preliminaries}\label{preliminaries}
\subsection{Flats of $\alpha$ and $\beta$}
Throughout this paper we shall pay special attention to the parts of the graph of $\alpha$ (resp. $\beta$) which are contained in proper affine subspaces of $H^1 (M,\R) \times \R$ (resp. $H_1(M,\R) \times \R$). Such sets are called flats of $\alpha$ (resp. $\beta$). We shall often identify a flat of $\alpha$ (resp. $\beta$) with its canonical projection to $H^1(M,\R)$ (resp. $H_1(M,\R)$). When $L$ is autonomous, by \cite{Carneiro}, $\alpha$ (not $\beta$) is constant on its flats. The sets of the kind $\partial \beta (h)$, $\partial \alpha (c)$ are the most obvious examples of flats.

We call relative interior of a flat of $\alpha$ (resp. $\beta$), its interior in the affine subspace it generates in $H^1 (M,\R) \times \R$ (resp. $H_1(M,\R) \times \R$). Recall Lemma A.3 of \cite{AvsM} : 
\begin{lemma}\label{nonor_4.2}
Let 
\begin{itemize}
         \item $E$ be a finite dimensional Banach space
	\item $A \co E \longrightarrow \R$ be a convex and superlinear map
	\item $x_0$ be a point of $E$ 
         \item $I$ be some (possibly infinite) set
         \item $F_i$, $i\in I$ be a family of flats of $A$ such that $x_0$ lies in the relative interior of $F_i$ for all  $i \in I$.
\end{itemize} 
Then there exists a flat $F$ containing $F_ i$ for all $i \in I$ such that  $x$ is an interior point of $F$.
\end{lemma}

This lemma enables us to speak of the largest flat of $\alpha$ containing a cohomology class $c$ in its relative interior. We denote it by  $F_c$.

\begin{lemma}\label{AvsM, A4}
Let 
\begin{itemize}
         \item $E$ be a finite dimensional Banach space
	\item $A \co E \longrightarrow \R$ be a convex and superlinear map
	\item $B  \co E^* \longrightarrow \R$ be the Fenchel dual of $A$.
\end{itemize}
Take  $x$ in  $E$ and  $y \in E^*$ in  $\partial A(x)$. Then any flat of $A$ containing $x$ in its interior is contained in $\partial B(y)$. In particular, if $x$ lies in the relative interior of $\partial B(y)$,  the largest flat of $A$ containing $x$ in its interior is   $\partial B(y)$.
\end{lemma}
For our purposes, this means that if $c$ is a cohomology class which lies in the relative interior of $\partial \beta (h)$, for some homology class $h$, then $F_c$, the largest flat of $\alpha$ containing $c$ in its relative interior,  is $\partial \beta (h)$.
\subsection{Aubry sets}\label{Aubry}
Define, for all $n \in \N$, 
	\[ 
	\begin{array}{rcl} h_n \co \left(M\times \T \right)\times \left(M\times \T \right) & \longrightarrow & \R \\
\left((x,t),(y,s)\right)& \longmapsto & \min \int^{s+n}_{t}L(\gamma,\dot\gamma,t)dt	
\end{array}
\]
where the minimum is taken over all absolutely continuous curves \\
$\gamma \co \left[t,s+n\right]\longrightarrow M$ such that $\gamma (t)=x$ and $\gamma (s+n)= y$. Note that this is a slight abuse of 
notation, since we denote by the same $t$ an element of $\T=\R/\Z$ or the corresponding point in $\left[0,1\right[$. 
The Peierls barrier is then defined as
	\[ \begin{array}{rcl} h \co \left(M\times \T \right)\times \left(M\times \T \right) & \longrightarrow & \R \\
\left((x,t),(y,s)\right)& \longmapsto & \liminf_{n \rightarrow \infty} h_n \left((x,t),(y,s)\right).
\end{array}
\]
The Aubry set is
	\[\tilde{\mathcal{A}}(L) := \left\{(x,t) \in M\times \T  \co h\left((x,t),(x,t)\right) =0 \right\}.
\]

 The image of the Aubry set under the canonical projection $(\pi, id) \co TM \times \T \longrightarrow M \times \T$, where $id$ is the identity map from $\T$ to itself,  is called projected Aubry set and denoted $\mathcal{A}(L)$. 
For the sake of brevity, we may sometimes say Aubry set rather than projected Aubry set. 
\subsubsection{Quotient Aubry sets}\label{Aubry quotient}
  Consider the equivalence relation on $\Azero$ defined by $(x,t)\approx (y,s)$ if and only if $h((x,t),(y,s))+h((y,s),(x,t))=0$. The classes are called static classes.The quotient $\Azero/\approx$ is a metric space with distance 
$$
d(\overline{(x,t)},\overline{(y,s)})=h((x,t),(y,s))+h((y,s),(x,t)),
$$ 
where $\overline{(x,t)}$ is the equivalence class of $(x,t)$. It is called  quotient Aubry set and denoted $A_0$ after \cite{Mather02}. The distance 
$d(\overline{(x,t)},\overline{(y,s)})$  is called the Mather distance. 

The  property of the quotient Aubry set that we use the most is the following  upper semi-continuity property : when the  quotient Aubry set of $L$ has 1-dimensional Hausdorff measure zero, by Corollary 5 of \cite{Bernard_Conley}, for any neighborhood $V$ of $\mathcal{A}(L)$ in $TM$, there exists a neighborhood $\mathcal{U}$ of $L$ in the $C^2$ compact-open topology such that for any $L'$ in $\mathcal{U}$, we have $\mathcal{A}(L') \subset V$. 

 The  quotient Aubry set of $L$ having 1-dimensional Hausdorff measure zero is not a rare phenomenon : by \cite{BC}, given a Tonelli Lagrangian $L$ on a closed manifold $M$, there exists a residual subset $\mathcal{O}(L)$ of $C^{\infty}( M)$, such that for every $f \in \mathcal{O}(L)$, for every $c \in H^1(M,\R)$, the quotient Aubry set of $(L,c)$ is a finite set. Also, by \cite{FFR}, for any autonomous Tonelli Lagrangian $L$ on a closed manifold of dimension two, the  quotient Aubry set of $L$ has 1-dimensional Hausdorff measure zero.

\subsection{Aubry sets and faces of $\alpha$}
As in the case of Mather sets, if $\omega$ is a closed one-form, $L-\omega$ has an Aubry set, which depends only on the cohomology class $c$ of 
$\omega$. We denote it $\mathcal{A}(L,c)$ or just $\mathcal{A}(c)$ when no confusion is possible. 
For the convenience of the reader we recall  Proposition 6 of \cite{ijm} (\cite{Bernard_02} for the time-periodic case)  :
\begin{proposition} \label{pomface}
Let $M$ be a closed manifold and let $L \co TM \times \T \longrightarrow \R$ be a  Tonelli Lagrangian. 
If a cohomology class $c_1$ belongs to 
a flat $F_{c}$ of $\alpha_L $ containing 
$c$ in its interior,
then $\tilde{\mathcal{A}}(c) \subset \tilde{\mathcal{A}}(c_1)$. 
In particular, if $c_1$ lies in 
the interior of $F_{c}$, then $\tilde{\mathcal{A}}(c) = \tilde{\mathcal{A}}(c_1)$. 
Conversely,
if two cohomology classes $c$ and $c_1$ are such that 
$\tilde{\mathcal{A}}(c) \cap \tilde{\mathcal {A}}(c_1) \neq \emptyset$, then $\alpha_L $ has a flat containing $c$ and $c_1$.
\end{proposition}
So for any flat $F$ of $\alpha$ and any $c_1$, $c_2$ in the relative  interior of $F$, the Aubry sets $\tilde{\mathcal{A}}(c_1)$ and 
$\tilde{\mathcal{A}}(c_2)$ coincide. We denote by $\tilde{\mathcal{A}}(F)$ the common Aubry set to all the cohomologies in the interior of $F$.  Recall that for any homology class $h$, $\partial \beta (h)$ is a flat of $\alpha$. For brevity we shall denote by $\tilde{\mathcal{A}}(h)$ the Aubry set 
$\tilde{\mathcal{A}}(\partial \beta (h))$, or $\tilde{\mathcal{A}}(L,h)$  when we need to emphasize the dependance on the Lagrangian.

\subsubsection{}\label{E_h et E_c}
 Now we look at the interplay between the dimension of $F_c$, and the size of the Aubry set of $c$. 

We call $\tilde{V}_c$ the  vector space  generated by  the pairs $(c'-c, \alpha (c')-\alpha (c))$ where $c' \in F_c$. 
We call $V_c$  the canonical projection of $\tilde{V}_c$ to $H^1 (M,\R)$. So, by Lemma \ref{AvsM, A4}, when $c$ lies in the relative interior of $\partial \beta (h)$ for some homology class $h$, we have $\tilde{V}_c= \tilde{V}_h$, with $\tilde{V}_h$ as  in Definition \ref{def V_h}.
\begin{definition} 
Let $\tilde{E}_c$ be the set of $(c',\tau) \in H^1 (M\times \T,\R)=H^1 (M,\R)\times H^1 ( \T,\R)$ such that there exists  a smooth closed one-form $\omega$ on $M\times \T$ with $[\omega] = (c',\tau)$ and $\spt (\omega)\cap \tilde{\mathcal{A}}(L,c) = \emptyset$. Let $E_c$ be the canonical projection of $\tilde{E}_c$ to $H^1 (M,\R)$.
\end{definition}
Note that if $c$ lies in the relative interior of $\partial \beta (h)$ for some homology class $h$, we have $\tilde{E}_c= \tilde{E}_h$, with $\tilde{E}_h$ as  defined in subsection  \ref{stronger stronger}.

The following theorem is proved in \cite{soussol}(see \cite{ijm} for the autonomous case). 
\begin{theorem}\label{ijm}  Let $L$ be a Tonelli Lagrangian on a closed manifold $M$, and let $c$ be a cohomology class in $H^1(M,\R)$. Then we have $E_c \subset V_c$.
\end{theorem}
Therefore, for every $h \in H_1 (M,\R)$, we have $E_h \subset V_h$, because every $E_h$ (resp. $V_h$) is an $E_c$ (resp. $V_c$) for some $c \in H^1(M,\R)$. 

Since the Differentiability Problem requires the integrality of $\tilde{V}_h$, rather than $V_h$, we first extend Theorem \ref{ijm} to $\tilde{E}_c$ and 
$\tilde{V}_c$. However we have to be careful because of the sign of $\alpha(c)$ in the respective definitions of $\tilde{E}_c$ and 
$\tilde{V}_c$.
Define a map 
$$
\begin{array}{rcl}
\mbox{flip} \co H^1(M,\R) \times \R & \longrightarrow  & H^1(M,\R) \times \R \\
(c,\tau) & \longmapsto & (c,-\tau).
\end{array}
$$
\begin{lemma}\label{flip}
 Let $L$ be a time-periodic Lagrangian on a closed manifold $M$, and let $c$ be a cohomology class in $H^1(M,\R)$. Then
 $$
 \tilde{E}_c \subset \mbox{flip} \left( \tilde{V}_c\right), 
 $$ 
and $E_c = V_c$ if and only if  $ \tilde{E}_c = \mbox{flip} \left( \tilde{V}_c\right)$.
 \end{lemma}
 \proof
 First we observe that $ \tilde{V}_c$ is the graph, over $V_c$, of the map $c' \mapsto \alpha (c')-\alpha (c)$. Now recall   Lemma 15 of \cite{soussol}:
\begin{lemma}\label{1forme}
If $\omega$ is a closed one form on $M\times \T$, with 
$[\omega]=(c,\tau) \in H^{1}(M,\R)\times H^{1}(\T,\R)$, and $\mu$ is an $(L,c)$-minimizing measure, then
\[
\int (L-\omega)d\mu = -\alpha (c) -\tau.
\]
\end{lemma}
Consequently, if $(c',\tau) \in \tilde{E}_c$, and $\omega$ is a smooth 1-form on $M\times \T$ with $[\omega] = (c',\tau)$ and $\spt (\omega)\cap \tilde{\mathcal{A}}(L,c) = \emptyset$, we have, for any $(L,c)$-minimizing measure $\mu$,
$$
-\alpha (c') -\tau=  \int (L-\omega)d\mu =  \int Ld\mu = -\alpha (c)
$$ 
so $\tau = \alpha (c)-\alpha (c')$. Therefore $ \tilde{E}_c$ is the graph, over $E_c$, of the map $c' \mapsto \alpha (c)-\alpha (c')$. Thus 

\begin{itemize}
  \item $E_c \subset V_c$ entails $\tilde{E}_c \subset \mbox{flip}\left( \tilde{V}_c\right)$, and since the former is always true, so is the latter. This proves the first statement of the lemma. 
  \item $E_c = V_c$ entails $\tilde{E}_c = \mbox{flip}\left( \tilde{V}_c\right)$, which proves the first implication of the second statement. The converse implication is obvious. 
\end{itemize}
 \qed

\begin{corollary}\label{E_c=V_c, V tilde entier}
 Let $L$ be a time-periodic Lagrangian on a closed manifold $M$, and let $c$ be a cohomology class in $H^1(M,\R)$. Then $E_c = V_c$ entails that $ \tilde{V}_c$ is  an integer subspace of $H^1(M,\R)\times \R$.
 \end{corollary}
 \proof
 By Lemma \ref{flip}, $E_c = V_c$ yields $ \tilde{E}_c = \mbox{flip} \left( \tilde{V}_c\right)$, and by Corollary \ref{E_0_entier}, $ \tilde{E}_c$  is  an integer subspace of $H^1(M,\R)\times \R$. Now the map $ \mbox{flip}$ leaves the integer lattice invariant, so it takes an integer subspace to an integer subspace. 
 
 \qed

In \cite{vmb2} we prove the following proposition:  
\begin{proposition} \label{prop vmb2}
Let $L$ be a Tonelli Lagrangian on a closed manifold $M$ 
with first Betti number $b$.
Assume that for every cohomology class $c$, $\tilde{V}_c$ is an integer subspace of 
$H^1 (M\times \T^1,\R)$. Let $h$ be a $k$-irrational homology class. 
Then $\beta$ is differentiable at $h$ in at least $k$ directions.
\end{proposition}
This implies that if $\tilde{V}_c$ is an integer subspace of 
$H^1 (M\times \T^1,\R)$ for every $c \in H^1 (M,\R)$, then we have an affirmative answer to the Differentiability Problem. However, in the proof of Proposition \ref{prop vmb2}, we only use cohomology classes which lie in the relative interior of $\partial \beta (h)$. So what we actually prove is 
\begin{proposition} \label{prop vmb2, modifiee}
Let $L$ be a Tonelli Lagrangian on a closed manifold $M$ 
with first Betti number $b$.
Assume that for every homology class $h$, $\tilde{V}_h$ is an integer subspace of 
$H^1 (M \times \T^1,\R)$. Let $h$ be a $k$-irrational homology class. 
Then $\beta$ is differentiable at $h$ in at least $k$ directions.
\end{proposition}
Proposition \ref{prop vmb2, modifiee} and Corollary \ref{E_c=V_c, V tilde entier} combine to prove  the following corollary, which  explains why the hypothesis  $E_h = V_h$ for any homology class $h$ in our main theorem contains an affirmative answer to the Differentiability Problem : 
\begin{corollary}\label{differentiabilite et E_h = V_h}
Let $M$ be a closed manifold and let $L \co TM\times \T$ be a Tonelli Lagrangian on $M$. Let $h \in H_1(M,\R)$ be a $k$-irrational homology class such that $E_h= V_h$. Then $\beta$ is differentiable at $h$ in at least $k$ directions.
\end{corollary}

\subsection{One last lemma before we go}\label{2.1}
Let us choose, for every $c \in H^1(M,\R)$, a closed one-form $\omega (c)$ on 
$M \times \T$,  such that  $\left[\omega\right]=(c,\alpha(c))$, in such a way that the map $c \mapsto \omega (c)$ is linear (hence continuous since $b_1(M)$ is finite). 

Recall that  $F_c$ is the maximal face of the epigraph of $\alpha$ that contains $(c,\alpha(c))$ in its relative interior. This definition makes sense by  Lemma \ref{nonor_4.2}. By Proposition \ref{pomface} we also have 
$$
 F_c (L) := \{ \left(c',\alpha(c')\right) \co \mathcal{A}(L,c) \subset \mathcal{A}(L,c') \}.
$$
  
 \begin{lemma}\label{C1}
 Let
\begin{itemize}
	\item $M$ be a closed manifold
	\item $L$ be a Tonelli Lagrangian on $TM \times \T$
	\item $c$ be a cohomology class in $H^1(M,\R)$.
\end{itemize}
Then 
	\[ F_0(L-\omega(c)) = F_c(L) -(c,0).
\]
\end{lemma}
\proof
First observe that for any cohomology class $c'$, 
\begin{eqnarray*}
\alpha_{L-\omega(c)}(c') & = & -\inf_{\mu} \int\left(L-\omega(c)-\omega(c') \right)d\mu \\
&=& -\inf_{\mu} \int\left(L-\omega(c+c') \right)d\mu \\
&=& \alpha_{L}(c+c').
\end{eqnarray*}
Now, for any cohomology class $c'$,
\begin{eqnarray*} 
 (c',\alpha_L(c')) & \in & F_c(L)  \Leftrightarrow  \\
\mathcal{A}(L,c) & \subset & \mathcal{A}(L,c') \Leftrightarrow  \\
\mathcal{A}(L-\omega(c)) & \subset & \mathcal{A}(L-\omega(c')) \Leftrightarrow \\
\mathcal{A}(L-\omega(c)) & \subset & \mathcal{A}(L-\omega(c)-\omega(c'-c)) \Leftrightarrow \\
(c'-c,\alpha_{L-\omega(c)}(c'-c)) & \in & F_0(L-\omega(c))  \Leftrightarrow  \\
(c'-c,\alpha_{L}(c')) & \in & F_0(L-\omega(c)) .
\end{eqnarray*}
\qed

\section{Lower semi-continuity  results}\label{semi-continuity}
When proving our main theorem we shall have to deal with the following situation : we have a sequence of homology classes $h_n$ that converges to some $h$. We know that the Aubry set of $h_n$ has some nice property $P$. So for any $c_n$ in the relative interior of $\partial \beta (h_n)$, the Aubry set of $c_n$ has property $P$. We would like to deduce that  any $c \in \partial \beta (h)$ lies in the closure of the set of cohomology classes whose Aubry sets have property $P$. For this we need to show that any $c \in \partial \beta (h)$ is a limit point of a sequence $c_n$ in in the relative interior of $\partial \beta (h_n)$. This is false in general but Lemma \ref{semi-con_2} below covers our needs. This is the reason why we include this rather technical section. 
\subsection{Lower semi-continuity of $F_c$}
The meaning of the next lemma is, roughly speaking, that under appropriate non-degenaracy hypothesis, the maximal face of $\alpha$ containing $c$ in its interior is lower semi-continuous as a function of $c$. 
\begin{lemma}\label{semi-con}
Let 
\begin{itemize}
\item $M$ be a closed manifold  
\item $c$ be a cohomology class in $H^1(M,\R)$
\item  $L_0$ be a Tonelli Lagrangian on $M$ 
 \item $\mathcal{A}_c$ be the Aubry set of $(L_0,c)$
\item $\alpha$ be the $\alpha$-function of $L_0$
\item $F_1$ be a compact, convex subset of the relative interior of $F_c (L_0)$
 containing $(c,\alpha(c))$ in its relative interior.
\end{itemize}
Assume that $E_c (L_0) = V_c (L_0)$,  and, for every cohomology class $c'$ in the relative interior of $F_c (L_0)$,  the quotient Aubry set of $(L_0,c')$ has Hausdorff one-dimensional measure zero.

Then there exists a neighborhood $\mathcal{U}$ of $L_0$ in the $C^2$ compact-open topology, such that for all $L \in \mathcal{U}$, denoting by 
$\alpha_L$ the $\alpha$-function of $L$, for all $c'$ such that $(c',\alpha(c')) \in F_1$, we have $(c',\alpha(c')-\alpha(c)+\alpha_L (c)) \in  F_c( L)$. 
\end{lemma}
\proof

Replacing $L_0$ with $L_0 -\omega$, where $\omega$ is any closed one-form with cohomology $c$, we  assume that $c=0$. For the sake of brevity we denote  $F_0 := F_0 (L_0)$.
Since $F_1$ is contained in the relative interior of 
$F_0$, by Proposition \ref{pomface}, for any $(c, \alpha(c)) \in F_1$, $\mathcal{A}(L_0,c)= \mathcal{A}(L_0)$. 
Since $E_0 (L)=V_0(L)$, there exists a neighborhood $U$ of $\mathcal{A}_0$ in 
$M \times \T$ such that for any  $(c,\alpha(c)) \in F_0$, there exists a closed one-form $\omega$ on 
$M \times \T$, supported outside of $U$, such that $\left[\omega\right]=(c,\alpha(c))$.

We have made the hypothesis  that for  any cohomology class $c$ in the relative interior of $F_0 (L_0)$, the one-dimensional Hausdorff measure of the quotient Aubry set of $L-\omega(c)$ is zero. So, by \cite{Bernard_Conley}, the Aubry set is semi-continuous, as a function of the Lagrangian, at $L_0 - \omega (c)$ for every $c$ in the relative interior of  $F_0$, where $\omega(c)$ is defined in Subsection \ref{2.1}. Thus, for every $c$ such that $(c,\alpha(c))$ lies  in the relative interior of  $F_0$, there exists a neighborhood $\mathcal{U}_1(c)$ of $L_0 $ in the $C^2$ compact-open topology, such that for all $L$ in $\mathcal{U}_1(c)$, we have $\mathcal{A}(L,c) \subset U$.

Now for any $(c, \alpha (c))$ in $F_1$, there exists a neighborhood  $V(c)$ of $c$ in $F_0$, and a neighborhood $\mathcal{U}(c)$ of $L_0 $ in the $C^2$ compact-open topology,  such that $L -\omega(c')+\omega(c)$ lies in $\mathcal{U}_1(c)$ for any $c'$ in $V(c)$ and any $L$ in $\mathcal{U}(c)$. This is where we use the fact that $F_1$ is contained in the relative interior of $F_0$.
Observe that for any $L,c,c'$, 
$$\mathcal{A}(L-\omega(c')+\omega(c),c)= \mathcal{A}(L,c'), $$
 so for all  $(c,\alpha(c)) \in F_1$, 
\begin{equation}\label{lemme 3.2, eq 1}
   \forall c' \in V(c), \  \forall L \in \mathcal{U}(c),\  
\mathcal{A}(L-\omega(c')+\omega(c),c)= \mathcal{A}(L,c') \subset U.
\end{equation}

Cover the compact set of $c$'s such that $(c,\alpha(c)) \in F_1$ by finitely many $V(c)$'s, say $V(c_1), \ldots V(c_n)$.
Then 
$$
\mathcal{U} := \bigcap_{i=1}^n  \mathcal{U}(c_i)
$$
is a neighborhood of $L_0$ in the $C^2$ compact-open topology. Take  any  
 $(c,\alpha(c)) \in F_1$. Let $i$ be such that $c \in V(c_i)$. Then  for any $L \in \mathcal{U}$, we have $L \in \mathcal{U}(c_i)$,
so by Equation (\ref{lemme 3.2, eq 1}) we have $\mathcal{A}(L,c)\subset U$. Recall that $E_0(L_0)$ is generated by 1-forms supported outside of $U$.
Thus for any $c$  such that $(c,\alpha(c)) \in F_1$,  for all $L \in \mathcal{U}$, 
$$
 E_0 (L_0) \subset E_c (L) \subset V_c (L).
$$
Recall that $V_c (L)$ is the underlying vector space of the affine space generated by $F_c (L)$,  
so for any $L \in \mathcal{U}$, for any $c$  such that $(c,\alpha(c)) \in F_1$, the graph of $\alpha_L$ contains an open subset of $(c,\alpha_L(c)) +E_0$. 

Now let us prove that for any $c$ such that $(c,\alpha(c))  \in F_1$, we have $\alpha_L (c) = \alpha (c) + \alpha_L(0) -\alpha (0)$. Pick any $c$ such that $(c,\alpha(c))  \in F_1$, and consider the map 
$$
\begin{array}{rcl}
\left[0,1\right] & \longrightarrow & \R \\
t & \longmapsto & \alpha_L (tc).
\end{array}
$$
This map has the same derivative as the map $t  \longmapsto \alpha (tc)$ because for any $t \in \left[0,1\right]$, the graph of $\alpha_L$ contains an open subset of $(tc,\alpha_L(tc)) +E_0$, while the  the graph of $\alpha$ contains an open subset of $(tc,\alpha (tc)) +E_0$. Thus, for all $t \in \left[0,1\right]$, we have $\alpha_L (tc) = \alpha (tc) + \alpha_L(0) -\alpha (0)$, in particular $\alpha_L (c) = \alpha (c) + \alpha_L(0) -\alpha (0)$.

Since this is true for any $c$ such that $(c,\alpha(c))  \in F_1$,  the graph of $\alpha_L$  contains $F_1 +(0,\alpha_L(0)-\alpha(0))$.
\qed
\subsection{Lower semi-continuity of $\partial \beta(L,h)$}
The next lemma says that $\partial \beta(L,h)$ is lower semi-continuous as a function of $L$ and  $h$, when $h$ is restricted to a special subspace of $H_1 (M,\R)$.
\begin{lemma}\label{semi-con_2}
Let
\begin{itemize}
\item $M$ be a closed manifold
\item $L$ be a Tonelli Lagrangian on $TM \times \T$
\item $\alpha$ denote the $\alpha$-function of $L$ 
\item $L_n$ a sequence of Lagrangians that converges to $L$ in the $C^2$ compact-open topology 
\item $h_0$ be a homology class in $H_1 (M,\R)$
\item $c_0$ be a cohomology class in the relative interior of  $\partial \beta(L,h_0)$
\item $H_0 := \{ h \in H_1(M,\R) \co \alpha(c)-\alpha(c') =\langle c-c',h \rangle \  \forall c,c' \in  \partial \beta(L,h_0) \}$
\item $h_n, n \in \N$ be a sequence in $H_0$  that converges to $h_0$
\item $c_n$ be an element of $\partial \beta(L_n, h_n)$ for each $n \in \N$, such that the sequence $c_n$ converges to some $c \in \partial \beta(L,h_0)$.
\end{itemize}
Assume that $E_{h_0}=V_{h_0}$, and, for every cohomology class $c$ in the relative interior of $\partial\beta (h_0)$,  the quotient Aubry set of $(L,c)$ has Hausdorff one-dimensional measure zero.

Then $c_n+c_0-c$ lies in $\partial \beta(L_n, h_n)$ for $n$ large enough.
\end{lemma}
\proof
We shall denote $\alpha_n$ and $\beta_n$ the $\alpha$ and $\beta$ functions of $L_n$, respectively.

First observe that 
\begin{eqnarray*}
\alpha (c_0) + \beta (h_0) & = & \langle c_0, h_0 \rangle \\
\alpha (c) + \beta (h_0) & = & \langle c, h_0 \rangle
\end{eqnarray*}
since $c,c_0 \in \partial \beta(h_0)$, so 
$\alpha(c_0)-\alpha(c)= \langle c_0 - c , h_0 \rangle  $. 
Besides, $c_n \in \partial \beta(L_n,h_n)$, so 
\begin{eqnarray*}
\alpha_n (c_n) + \beta_n (h_n) & = & \langle c_n, h_n \rangle \\
 & = & \langle c_n +c_0 -c, h_n \rangle + \langle c - c_0 , h_n \rangle \\ 
 &\leq & \alpha_n (c_n +c_0 -c) + \beta_n (h_n)  +\alpha(c)-\alpha(c_0)
\end{eqnarray*}
where we have used the Fenchel inequality for $c_n +c_0 -c$ and $h_n$, and the fact that $h_n \in H_0$. Therefore
\begin{equation}\label{eq1}
\alpha_n (c_n) +\alpha(c_0)-\alpha(c) \leq  \alpha_n (c_n +c_0 -c) .
\end{equation}
We shall now prove, for $n$ large enough,  the converse inequality
\begin{equation}\label{eq2}
\alpha_n (c_n) +\alpha(c_0)-\alpha(c) \geq  \alpha_n (c_n +c_0 -c).
\end{equation} 
It will follow that 
 \begin{eqnarray*}
\alpha_n (c_n+c_0-c) + \beta_n (h_n) &=& \alpha_n (c_n) +\alpha(c_0)-\alpha(c)+ \beta_n (h_n)\\
&=& \langle c_n , h_n \rangle +   \langle c - c_0 , h_n \rangle  \\
&=& \langle c_n+c_0-c , h_n \rangle \\
 \end{eqnarray*}
 using the Fenchel equality for $c_n$ and $h_n$, and the fact that $h_n \in H_0$. Therefore  $c_n+c_0-c \in \partial \beta(L_n,h_n)$, which proves the lemma. 
 
 Now let us  prove Equation (\ref{eq2}).
  Since $c_0$ lies in the relative interior of $\partial \beta(h_0)$, by  Lemma \ref{AvsM, A4}, 
  \[
  F_{c_0}(L)= \left\{(c, \alpha_L (c)) \co c \in \partial \beta(h_0) \right\}.
  \] 
  Since  $c$ lies in $\partial \beta(h_0)$ and $c_0$ lies in the relative interior of $\partial \beta(h_0)$, there exists a positive $\epsilon$ such that 
$$
\{ \left( c_0 +t(c-c_0), \alpha(c_0) + t(\alpha(c)-\alpha(c_0)) \right) \co t \in \left[-2\epsilon, 1 \right] \} \subset F_{c_0}(L).
$$
Therefore
$$
F_1 := \{ (c_0,\alpha (c_0)) + t \left( c-c_0, \alpha(c)-\alpha(c_0) \right) \co t \in \left[-\epsilon, \frac{1}{2} \right]\}
$$
is contained in the relative interior of $F_{c_0}(L)$. So  by Lemma \ref{semi-con}, there exists $N \in \N$ such that for all $n \geq N$, 
\begin{equation} \label{eq5}
F_1 + (0,\alpha_{L_n  -\omega(c_n) +\omega(c)} (c_0)-\alpha(c_0))\subset F_{c_0}(L_n-\omega(c_n) +\omega(c))
\end{equation}
where $\omega(c_n), \omega(c)$ are defined as in \ref{2.1}.
Thus, taking $t=1/2$ and $t=0$, we have, $\forall n \geq N$,
\begin{eqnarray*}
 \ \alpha_{L_n  -\omega(c_n) +\omega(c)}(\frac{c_0+c}{2})-\alpha_{L_n  -\omega(c_n) +\omega(c)}(c_0)
 &=& \alpha (\frac{c_0 + c}{2}) -\alpha (c_0)\\
 &=& \frac{\alpha(c) -\alpha(c_0)}{2}  
\end{eqnarray*}
and, recalling that $\alpha_{L_n  -\omega(c_n) +\omega(c)}(.) = \alpha_n(.+c_n-c)$,
\begin{equation}\label{4}
\forall n \geq N, \  \alpha_n(c_n + \frac{c_0 -c}{2}) -\alpha_n(c_n+c_0 -c) = 
\frac{1}{2}\left( \alpha(c)- \alpha(c_0)\right).
\end{equation}

 Combining Equation (\ref{eq5}) 
 and the fact that 
 \[
 \left(c_0, \alpha(c_0)\right), \left(\frac{c_0+c}{2}, \frac{\alpha(c_0)+\alpha(c)}{2}\right) \in F_1,
 \]
 we get that 
\begin{eqnarray*}
  \left(c_0, \alpha(c_0)\right) + (0,\alpha_{L_n  -\omega(c_n) +\omega(c)} (c_0)-\alpha(c_0)) &=&  
   \left(c_0, \alpha_{L_n  -\omega(c_n) +\omega(c)} (c_0)\right) \\
   &=&  \left(c_0, \alpha_{n  } (c_0+c_n-c)\right) 
 \end{eqnarray*}
   and 
   \begin{eqnarray*}
  && \left(\frac{c_0+c}{2}, \frac{\alpha(c_0)+\alpha(c)}{2}\right)  + (0,\alpha_{L_n  -\omega(c_n) +\omega(c)} (c_0)-\alpha(c_0))  \\
  &&
  =\left(\frac{c_0+c}{2}, \alpha_{L_n  -\omega(c_n) +\omega(c)} (c_0)+\frac{\alpha(c)-\alpha(c_0)}{2}\right) \\
  &&
 =\left(\frac{c_0+c}{2}, \alpha_{n } (c_0+c_n-c)+\frac{\alpha(c)-\alpha(c_0)}{2}\right)    
  \end{eqnarray*}
 both lie in   $ F_{c_0}(L_n-\omega(c_n) +\omega(c)) $.  
 Recall that by Lemma \ref{C1} 
	\[F_{c_n+c_0-c}(L_n)-(c_n -c,0)= F_{c_0}\left(L_n-\omega(c_n)+\omega(c)\right).
\]   
Therefore $  \left(c_0+c_n-c, \alpha_{n  } (c_0+c_n-c)\right) $
   and 
   $$
  \left(c_n + \frac{c_0-c}{2}, \alpha_{n } (c_0+c_n-c)+\frac{\alpha(c)-\alpha(c_0)}{2}\right)  
  $$   
  both lie in   $F_{c_n+c_0-c}(L_n)$. 

Now  let us  take some $h'_n$ in $\partial \alpha(L_n,c_n+c_0-c)$. First let us observe that for any $c'$ such that $(c',\alpha_n (c')) \in F_{c_n+c_0-c}(L_n)$, we have $c' \in  \partial \beta(L_n,h'_n) $. Indeed the projection to $H^1(M,\R)$ of $ F_{c_n+c_0-c}(L_n)$ is a flat of $\alpha$ containing $c_n+c_0-c$ in its relative interior, and $h'_n$ lies in $\partial \alpha(L_n,c_n+c_0-c)$.
Then, applying Lemma \ref{AvsM, A4}, with $x=c_n+c_0-c$ and $y=h'_n$, we get that the projection to $H^1(M,\R)$ of $F_{c_n+c_0-c}(L_n)$ is contained in $ \partial \beta(L_n,h'_n) $, which proves our claim.

Hence   both $c_n+c_0-c$ and  $c_n +2\inv (c_0 -c)$ lie in  
$\partial \beta(L_n,h'_n)$, that is, 
\begin{eqnarray*}
\alpha_n(c_n+c_0-c) + \beta_n (h'_n) &=& \langle c_n+c_0-c, h'_n \rangle \\
 \alpha_n \left(c_n + \frac{c_0 -c}{2}\right) + \beta_n (h'_n) &=& \langle c_n + \frac{c_0 -c}{2}, h'_n \rangle
 \end{eqnarray*}
 whence
$$
\langle \frac{c_0 -c}{2}, h'_n \rangle =
\alpha_n(c_n+c_0-c)- \alpha_n \left(c_n + \frac{c_0 -c}{2}\right) 
$$
and, using Equation (\ref{4}),  $\forall n \geq N, \   \langle c_0 -c, h'_n \rangle =  \alpha(c_0)- \alpha(c)$ (that is, $h'_n \in H_0$).
Therefore
\begin{eqnarray*}
\alpha_n(c_n) +\beta_n(h'_n) & \geq & \langle c_n , h'_n \rangle \\
&=& \langle c_n +c_0 -c, h'_n \rangle  +  \langle c -c_0, h'_n \rangle \\
&=& \alpha_n(c_n+c_0 -c) +\beta_n(h'_n) +\alpha(c)- \alpha(c_0)
\end{eqnarray*}
which proves Equation \ref{eq2}, and the lemma.
\qed

\section{Density of Legendre transforms of rational homologies}\label{density}
Theorem \ref{main} is an immediate corollary of the following 
\begin{theorem}\label{rationnel_dense_periodique}
Let 
\begin{itemize}
	\item $M$ be a closed manifold
	\item $L$ be a Tonelli Lagrangian on $TM \times \T$ 
	\item $U$ be an open set of $H^1(M,\R)$, such that for all $h$ in the Legendre transform $V := \partial \alpha(U)$ of $U$, $E_h =V_h$ and the quotient Aubry set $A_h$ has Hausdorff one-dimensional measure zero.
\end{itemize}
 Then the Legendre transform $V = \partial \alpha(U)$ contains a rational homology class.
\end{theorem}
\proof
Recall that $V$, the Legendre transform of $U$, is the set of homology classes $h$ such that for some $c \in U$, $<c,h>=\alpha(c)+\beta(h)$.
In Appendix \ref{paragraphe rationnel affine} we define a rational affine subspace subspace of $H_1 (M,\R)$ as a subset of $H_1 (M,\R)$ defined by affine equations with integer coefficients. 
 We shall  prove by induction on $k=0,1,\ldots b_1(M)-1$ the following alternative :
\begin{itemize}
	\item either $V$ contains an open subset of a rational affine subspace $H_k$ of $H_1 (M,\R)$, of codimension $k$
	\item or there exists an open subset $U_k$ of $U$, and integer one-forms  $\omega_1,\ldots \omega_{k+1}$ on $M\times \T$, whose cohomology classes are linearly independant in $H^1(M\times \T, \R)$, and  such that 
	\[ \forall c \in U_k,\  \forall i=1,\ldots k+1,\  \mathcal{A}_c \cap  \mbox{supp}\omega_i = \emptyset.
\]
\end{itemize}
In the first case, by Lemma \ref{lemme rationnel affine},  $V$ contains a rational homology class, so we are done. Assume we are in the second  case for $k=b_1(M)-1$. Set 
$\left[\omega_i\right] = (c_i, \tau_i) \in H^1(M,\R) \times H^1(\T,\R)$ for each 
$i=0,\ldots b_1(M)-1$.  Pick $c \in U_k$ and an $(L,c)$-minimizing measure $\mu$. We have $\int \omega_i d\mu =0$ because 
$\spt \mu \subset \mathcal{A}(c)$ and $ (\spt \omega_i ) \cap \mathcal{A}(c)=\emptyset$.
On the other hand
$
\int \omega_i d\mu = \langle c_i, h \rangle + \tau_i.
$
Now $\tau_i \in \Z$ because $\omega_i$ is an integer one-form, so $\langle c_i, h \rangle \in \Z$ for $i=1, \ldots b_1(M)$.

 Note that $\langle c_i, h \rangle $,  $0=1, \ldots b_1(M)-1$, are the coordinates of $h$ in the basis of $H_1(M,\R)$ dual to the  basis $c_i, i=0, \ldots b_1(M)-1$ of $H^1(M,\R)$. This basis  consists of integer cohomology classes, so its dual consists of rational homology classes. This proves that $h$ is rational, and the proposition.

Let us start the induction with $k=0$. 

\textbf{First case.} Assume that for all $h$ in $V$, $\dim  \partial \beta(h)=0$, that is, 
$\partial \beta(h)$ is a point. Then let us show that $V$ is open in $H_1 (M,\R)$. Take
\begin{itemize}
	\item $h \in V$
	\item a sequence $h_n$ in $H_1 (M,\R)$ such that $h_n$ converges to $h$
	\item a sequence $c_n$ in $H^1 (M,\R)$ such that $<c_n,h_n>=\alpha(c_n)+\beta(h_n)$ for all $n \in \N$.
\end{itemize}
Since $\alpha$ is superlinear, the sequence $c_n$ remains within some compact subset of $H^1 (M,\R)$, hence we may assume that $c_n$ converges to some $c$. Then by continuity $<c,h>=\alpha(c)+\beta(h)$ so $c \in \partial \beta(h)$. Since we assumed that $\partial \beta(h)$ is a point, we have $\partial \beta(h)= \left\{c\right\}$. Now recall that $h \in V$, so $c \in U$. Hence $\exists n_0 \in \N,\  \forall n \geq n_0,\  c_n \in U$. Therefore $\forall n \geq n_0,\  h_n \in V$, which proves that $V$ is open in $H_1 (M,\R)$, hence it contains a rational homology class.

\textbf{Second case.} Assume that for some $h$ in $V$, $\dim \partial \beta(h) \geq 1$. Since $h \in V$, we have $\partial \beta(h)  \cap U \neq \emptyset$, and since $U$ is open, $U$ must then meet the relative interior of $\partial \beta(h)$. Take $c \in U$ in the relative interior of $\partial \beta(h)$. Then the Aubry set, (resp.  quotient Aubry set), of $c$, are the Aubry set  (resp.  quotient Aubry set) of $h$. In particular the quotient Aubry set $A_c$ has Hausdorff one-dimensional measure zero.

Furthermore, $F_c$, the largest face of $\alpha$ containing $c$ in its relative interior, is $\partial \beta(h)$ by  Lemma \ref{AvsM, A4}. Thus $V_c = \mbox{Vect}\partial \beta(h)$, which yields $E_c = \mbox{Vect}\partial \beta(h)$ by  our assumption on $U$. So the dimension of $E_c$ is at least one. Moreover $E_c$ is an integer subspace of $H^1(M,\R)$ by Lemma \ref{E_0_entier}. So we may find an integer one-form $\omega_1$ on $M \times \T$,  such that $\left[\omega_1\right] \in E_c$ and the support of $\omega_1$ is disjoint from the Aubry set of $c$. 

Moreover, by the upper semi-continuity of the Aubry set, there exists a neighborhood $U_1$ of $c$ in $U$, such that 
$$
\forall c' \in U_1, \  \mathcal{A}(c') \cap  \mbox{supp} \omega_1 = \emptyset.
$$

 This finishes the first induction step.

Assume now we have carried out the induction process until the $k$-th step for some $1\leq k \leq b_1(M)-2$. If we are in the first case of the $k$-th step, there is nothing left to do, so assume we are in the second case. 

Set $\left[\omega_i\right] = (c_i, \tau_i) \in H^1(M,\R) \times H^1(\T,\R)$ for each 
$i=0,\ldots k$. Let $H_k$ be the rational affine  subspace of $H_1 (M,\R)$ 
defined by the equations $\langle c_i, . \rangle = -\tau_i $ for $i=1,\ldots k+1$. Denote by $V_k$  the Legendre transform of $U_k$. Then any $h \in V_k$ is the homology class of a $c$-minimizing measure $\mu$ for some $c \in U_k$. The support of $\mu$ does not meet $\mbox{supp} \omega_i$ so $\langle c_i, \left[\mu \right] \rangle = -\tau_i $ for $i=1,\ldots k+1$. 
Moreover, 
$$
\forall h \in V_k, \ \mbox{Vect}(c_1, \ldots c_{k+1}) \subset E_c \subset V_c \subset \mbox{Vect}\partial \beta(h) ,
$$
 where the last inclusion holds because $c \in \partial \beta (h)$, whence $F_c \subset  \partial \beta (h)$ by Lemma \ref{AvsM, A4}.
 Thus the dimension of $\partial \beta(h) $ is $\geq k+1$.

\textbf{First case.} Assume that for all $h \in V_k$, $\dim \partial \beta(h) = k+1$, that is, for all $h \in V_k$, $\mbox{Vect}(c_1, \ldots c_{k+1}) = \mbox{Vect}\partial \beta(h) $. Let us show, then, that $V_k$ is open in $H_k$. Take
\begin{itemize}
	\item $h_0 \in V_k$
	\item $c_0 \in U_k$ such that $c_0$ lies in the relative interior of $\partial \beta(h_0)$
	\item a sequence $h_n$ in $H_k$ such that $h_n$ converges to $h_0$
	\item a sequence $c_n$ such that  $c_n \in \partial \beta(h_n)$ for all $n\in \N$
	\end{itemize}
Taking a subsequence if we have to, we may assume the sequence $c_n$ converges to some $c$ in $\partial \beta(h_0)$. We want to apply Lemma \ref{semi-con_2} to $h_n$ so we have to check that  for every $c,c'$ in $\partial \beta (h_0)$, we have 
\begin{equation}\label{h_n dans H_0}
\alpha(c) -\alpha (c') = \langle c-c', h_n \rangle.
\end{equation}
Take $c,c'$ in $\partial \beta (h_0)$, so we have $\alpha(c) -\alpha (c') = \langle c-c', h_0 \rangle$. Since $\mbox{Vect}(c_1, \ldots c_{k+1}) = \mbox{Vect}\partial \beta(h_0) $, there exist real numbers $\lambda_1, \ldots \lambda_{k+1}$ such that $c-c' = \sum_{i=1}^{k+1} \lambda_i c_i$, so 
$ \langle c-c', h_0 \rangle =  -  \sum_{i=1}^{k+1} \lambda_i \tau_i$. Now since $h_n \in H_k$, we have, for all $i=1, \ldots, k+1$, $\langle c_i, h_n \rangle = -\tau_i$, so $ \langle c-c', h_n \rangle =  -  \sum_{i=1}^{k+1} \lambda_i \tau_i = \alpha(c) -\alpha (c')$, which proves Equation (\ref{h_n dans H_0}).

Then by Lemma \ref{semi-con_2} there exists  $N \in \N$ such that $\forall n \geq N$, $c_n + c_0 -c \in \partial \beta(h_n)$. Now $c_n + c_0 -c$ converges to $c_0$ so for $n$ large enough, $c_n + c_0 -c \in U_k$. Then $h_n \in V_k$, which proves that $V_k$ is open in $H_k$.

 \textbf{Second case.} Assume that for some $h \in V_k$, $\dim \partial \beta(h) > k+1$. Take $c$ in the relative interior of $\partial \beta(h)$. We have $E_c = \mbox{Vect}\partial \beta(h)$ as in the second case of the first step so the dimension of $E_c$ is at least $k+2$. Since $E_c$ is an integer subspace of $H^1(M,\R)$, we may find linearly independant integer one-forms $\omega_1, \ldots \omega_{k+2} $ such that 
 $$
\forall i= 1, \ldots k+2, \    \mathcal{A}(c) \cap  \mbox{supp} \omega_i = \emptyset.
$$  
 Moreover, by semi-continuity of the Aubry set, there exists a neighborhood $U_{k+1}$ of $c$ in $U_k$, such that 
$$
\forall c' \in U_{k+1}, \  \forall i= 1, \ldots k+2, \    \mathcal{A}(c') \cap  \mbox{supp} \omega_i = \emptyset.
$$  
 This finishes the $(k+1)$-th induction step, and the proof of the theorem. 
\qed

\appendix
\section{Integrality of $E_0$}
The following lemma is not directly useful for the proof of Theorem \ref{main}, only Lemmata \ref{E_0_entier} and \ref{lemme rationnel affine} are. However we believe it might be useful in future work on the subject, so we beg the reader to bear with us for a while. 
\begin{lemma}\label{topologie}
Let \begin{itemize}
  \item $M$ be a closed manifold, equipped with a Riemann metric
  \item $b_1(M)$ be the first Betti number of $M$
  \item $F$ be a closed subset of $M$
  \item $F_{\epsilon}$ be the open $\epsilon$-neighborhood of $F$ in $M$, for any positive $\epsilon$
  \item $E_{\epsilon}$ be the set of cohomology classes of  closed 1-forms on $M$ supported outside $F_{\epsilon}$
  \item $E_0$ be  the  set of cohomology classes of  closed 1-forms on $M$ supported outside $F$
  \item $H_{\epsilon}$ be the subset of $H_1 (M,\R)$ that consists of the homology classes of   cycles contained in $F_{\epsilon}$
  \item  $E_0^{\perp}$ be the set of homology classes $h$ in $H_1 (M,\R)$  such that $\langle c,h\rangle =0$ for all $c$ in $E_0$.
      \end{itemize}
     Then there exists $\epsilon_0 >0$ such that  for any $\ 0 < \epsilon \leq \epsilon_0$, any element of   $E_0^{\perp}$ is represented by a  cycle contained in        $F_{\epsilon}$, with at most $b_1(M)$ connected components.
     \end{lemma}
     \proof
Observe that 
$$
\forall \ 0 < \epsilon' \leq \epsilon,\  E_{\epsilon} \subset E_{\epsilon'} \mbox{ and } H_{\epsilon'} \subset H_{\epsilon}.
$$
Also, since $E_{\epsilon}$ (resp. $H_{\epsilon}$) is a vector subspace of $H^1 (M, \R)$ (resp. $H_1 (M, \R)$), which is finite dimensional since $M$ is compact, there exists $\epsilon_0 >0$ such that 
$$
\forall \ 0 < \epsilon \leq \epsilon_0 ,\   E_{\epsilon} = E_{\epsilon_0} \mbox{ and } H_{\epsilon} = H_{\epsilon_0}.
$$
Since $F$ is compact, any element of $E_0$ is contained in $E_{\epsilon}$ for some $\epsilon >0 $, so $E_0 = E_{\epsilon_0}$. 

Denote by $H_{\epsilon}^{\perp}$ the set of cohomology classes $c$ in $H^1(M, \R)$ such that $\langle c,h\rangle =0$ for all $h$ in $H_{\epsilon}$.

First let us show that $H_{\epsilon}^{\perp}= E_0$ for any $\ 0 < \epsilon \leq \epsilon_0$.
Take \begin{itemize}
  \item $\ 0 < \epsilon' < \epsilon \leq \epsilon_0$
   \item an element $c$ of   $H_{\epsilon}^{\perp}$   
   \item  a  closed 1-form $\omega$ on $M$ such that $\left[\omega \right]=c$
   \item a smooth function $\varphi$ on $M$ such that $\varphi (x)=1$ for all $x$ in $F_{\epsilon'}$, and  $\varphi (x)=0$ for all $x$ in $M \setminus    F_{\epsilon}$.   \end{itemize}
Then the integral of $\omega$ vanishes on any cycle contained in $F_{\epsilon}$, hence 
$\omega$ is exact inside $F_{\epsilon}$, that is, there exists a $C^1$ function
 $f \co F_{\epsilon} \longrightarrow \R$ such that $\omega = df$ inside $F_{\epsilon}$.

So the closed 1-form $\omega -d(\varphi f)$ is cohomologous to $\omega$, and vanishes identically inside $F_{\epsilon'}$, so $c \in E_{\epsilon'} =E_0$.
Therefore $H_{\epsilon}^{\perp} \subset E_0$. The converse inclusion is obvious, so 
$$
\forall \  0  < \epsilon \leq \epsilon_0, \   H_{\epsilon}^{\perp} = E_0  .
$$
Hence  by duality (recall that the dimension of $H_1(M,\R)$ is finite)
$$
\forall \  0  < \epsilon \leq \epsilon_0, \   H_{\epsilon} = E_0^{\perp},
$$
that is, any element $h$ of $H_1 (M, \R)$ such that $\langle c,h \rangle = 0$, $\forall c \in E_0$,  is represented by a cycle contained in $F_{\epsilon}$ for any $ 0  < \epsilon \leq \epsilon_0$.

We still have to prove the statement about the number of connected components. Consider the map $J \co H_1(F_{\epsilon}, \R) \longrightarrow  H_{\epsilon}$ induced by the inclusion of $F_{\epsilon}$ into $M$. The map $J$ is surjective by definition of $H_{\epsilon}$. The connected cycles contained in $F_{\epsilon}$ generate $H_1(F_{\epsilon},\R)$, hence they generate $H_{\epsilon}$. Therefore we may find a basis of $H_{\epsilon}$ that consists of connected cycles contained in $F_{\epsilon}$. The cardinal of this basis is at most $b_1(M)$ since $H_{\epsilon}$ is a vector subspace of $H_1 (M, \R)$. Therefore any element of $H_{\epsilon}$ is represented by a linear combination of at most $b_1(M)$ connected cycles contained in $F_{\epsilon}$. 

\qed

\begin{lemma}\label{E_0_entier}
Let
\begin{itemize}
 \item $M$ be a compact manifold without boundary
   \item $F$ be a closed subset of $M$ 
  \item $E_0$ be  the  set of cohomology classes of  closed one forms on $M$ supported outside $F$.
  \end{itemize}
  Then $E_0$ is an integer subspace of  $H^1(M,\R)$.
  \end{lemma}
  \proof
  Re-using the notation of the previous lemma, we have $H_{\epsilon}^{\perp}= E_0$ for any $\ 0 < \epsilon \leq \epsilon_0$. Now by the Universal Coefficient Theorem, $H_1(F{\epsilon}, \R)$ is generated by integer classes, so $H_{\epsilon}$ is an integer   subspace of $H_1(M,\R)$.  Therefore $H_{\epsilon}^{\perp}= E_0$ is also integer.
  \qed  
  
\subsection{Rational affine subspaces}\label{paragraphe rationnel affine}
We say an affine subspace of $H_1(M,\R)$ is rational if it is defined by equations of the form $\left\langle c_i, h \right\rangle =  \tau_i$, $i=1, \ldots k$, where  $c_i$, $i=1, \ldots k$, are integer cohomology classes, and $\tau_i \in \Z$, $i=1, \ldots k$. What we need to know about rational affine subspaces is the
\begin{lemma}\label{lemme rationnel affine}
Let  $H$ be a rational affine subspace of $H_1(M,\R)$.Then $H \cap H_1 (M,\Q)$ is dense in $H$.
\end{lemma}
\proof
Let $\left\langle c_i, h \right\rangle =  \tau_i$, $i=1, \ldots k$ be the equations that define $H$. Discarding some equations if we have to, we may assume that $c_i$, $i=1, \ldots k$ are linearly independant. Take integer cohomology classes $c_{k+1},  \ldots c_b$ such that $c_1,\ldots c_b$ is a basis of $H^1(M,\R)$ as a vector space. Then the numbers $\left\langle c_i, h \right\rangle$ are the coordinates of $h$ in the basis of $H_1(M,\R)$ dual to $c_1,\ldots c_b$, which consists of rational homology classes. Then the homology classes $h$ which satisfy
\begin{eqnarray*}
\left\langle c_i, h \right\rangle &=&  \tau_i,  \  i=1, \ldots k \\
\left\langle c_i, h \right\rangle & \in &  \Q,  \  i=k+1, \ldots b
	\end{eqnarray*}
are rational, and they form a dense subset of $H$.
\qed

  {\small

\bigskip

\noindent

D\'epartement de Math\'ematiques, Universit\'e Montpellier 2, France\\
e-mail : massart@math.univ-montp2.fr
}

\end{document}